\documentclass[a4paper,12pt]{article}

\usepackage{amssymb,amsfonts,amsmath,epsfig,float,graphicx}
\usepackage{hyperref} 

\newtheorem{theorem}{Theorem}[section]

\newtheorem{proposition}[theorem]{Proposition}
\newtheorem{definition}[theorem]{Definition}

\newcommand{\qed}{\nobreak \ifvmode \relax \else
      \ifdim\lastskip<1.5em \hskip-\lastskip
      \hskip1.5em plus0em minus0.5em \fi \nobreak
      \vrule height0.75em width0.5em depth0.25em\fi}

\def\upd{\,{\rm d}}
\newcommand{\be}{\begin{eqnarray}}
\newcommand{\ee}{\end{eqnarray}}
\newcommand{\nnee}{\nonumber\ee}
\newcommand{\Tr}{\,{\rm Tr}\,}
\newcommand{\dom}{\,{\rm dom}\,}

\def\Co{{\mathbb C}}

\def\Io{{\mathbb I}}

\def\Mo{{\mathbb M}}

\def\Ro{{\mathbb R}}

\def\beginproof{\par\strut\vskip 0.0cm\noindent{\bf Proof}\par}
\def\endproof{\par\strut\hfill$\square$\par\vskip 0.2cm}

\title{
Quantum statistical manifolds: the linear growth case
}
\author{
Jan Naudts\\
\strut\\
\small          Departement Fysica, Universiteit Antwerpen,\\
\small          Universiteitsplein 1, 2610 Antwerpen, Belgium\\
\small          \url{Jan.Naudts@uantwerpen.be}\\
\small          \url{https://orcid.org/0000-0002-4646-1190}
}

\date {}

\begin{document}
\maketitle

\begin{abstract}
A class of vector states on a von Neumann algebra is constructed.
These states belong to a deformed exponential family.
One specific deformation is considered.
It makes the exponential function asymptotically linear.
Difficulties arising due to non-commutativity are highlighted.
\end{abstract}

\section{Introduction}

In a recent publication Montrucchio and Pistone \cite{MP17} treat a special case of a 
para\-meter-free deformed
exponential family of probability distributions.
The present paper shows that part of this work
can be transposed to a non-commutative setting in a rather straightforward manner.
Both the commutative and the non-commutative versions can be useful as an inspiration
for the development of a more general theory of para\-meter-free information geometry.
It is not the ambition of the present paper to develop such a theory, but only to
clarify the kind of difficulties which one encounters in the non-commutative setting.

Amari \cite{AS85,AN00} studied parametrized models of Information Geometry.
Para\-meter-free families were introduced by Pistone and Sempi \cite{PS95}.
See also \cite {ZH06,AJVS17}.
The generalization of Information Geometry to a non-commutative context is of interest because of
its applications in Quantum Theory. However, it is not so trivial.
In the simplest context the random variables of probability theory
are replaced by $n\times n$-dimensional matrices,
the probability distribution is replaced by a density matrix. Then use
can be made of the property known as cyclic permutation under the trace.
This property restores part of the commutativity, needed to mimic the proofs
of the commutative case. A more general context involves Tomita-Takesaki theory.
See for instance the recent book of Petz \cite{PD08}. The aim of the present paper
is to go beyond the traditional setting by not longer focusing on tracial states.
In the Gelfand-Naimark-Segal (GNS) representation one can make use
of the commutant algebra. This relaxes problems with non-commutativity as well.
A treatment of the matrix case along these lines  has been tried out in \cite{NJ18}.
A study of log-affine geodesics in a manifold of states on a von Neumann algebra 
is found in \cite {NJ19}.

A regularization of the exponential function was introduced by Newton \cite{NNJ12}.
The idea was picked up by Montrucchio and Pistone \cite{MP17}.
The deformation of the exponential function is used to construct deformed exponential families
of probability distributions and, in the present paper, of quantum states.

The interest in deformed exponential families started with the q-statistics of
Tsallis \cite{TC88}. A further generalization was given by the author \cite{NJ02,NJ04,NJ05,NJ11}.
The latter formalism is used here and is explained below in Section \ref {sect:deffun}
for the special case of linear growth.

\subsection*{Non-commutative context}

A statistical manifold is a differentiable manifold  $\Mo$
together with a Riemannian metric $g$ and a pair of dually flat
connections \cite{AN00}. The manifold $\Mo$ consists of
probability distributions on a given measure space $({\cal X},\upd x)$.

In the most simple non-commutative setting the elements of $\Mo$ are
density matrices instead of probability distributions.
These are self-adjoint matrices with non-negative eigenvalues and with trace 1.
In the present work the more general $C^*$-algebraic context is chosen.
A state $\omega$ on a $C^*$-algebra $\cal A$ is a linear map $A\in {\cal A}\mapsto \omega(A)\in\Co$
which satisfies the positivity condition
that $A\ge 0$ implies $\omega(A)\ge 0$ and the normalization
condition $\omega(\Io)=1$ (for convenience, it is assumed that the identity $\Io$ belongs to $\cal A$).
Note that any density matrix $\rho$ determines a state $\omega$ of the $C^*$-algebra
of all square matrices $A$ of given dimension by the relation $\omega(A)=\Tr \rho A$.

Given a state $\omega$ on the $C^*$-algebra ${\cal A}$ there exists a *-representation $\pi$
of ${\cal A}$ as bounded linear operators on a Hilbert space ${\cal H}$, together with an
element $\Omega$ of ${\cal H}$ such that 
\be
\omega(A)&=&(\pi(A)\Omega,\Omega)\quad \mbox{for all}\quad A\in {\cal A}
\ee
and such that $\pi({\cal A})$ is dense in ${\cal H}$. 
This representation is unique up to unitary
equivalence. Its is known as the Gelfand-Naimark-Segal (GNS)
representation induced by the state $\omega$.

Let us make the simplifying assumptions that the $C^*$-algebra $\cal A$ is a von Neuman algebra
of operators on a fixed Hilbert space $\cal H$ and 
that there is given a fixed faithful normal state $\omega$,
which will be used as starting point of the construction following later on.
Because $\omega$ is faithful there exists
an element $\Omega$ of $\cal H$ such that the trivial representation,
defined by $\pi(A)=A$ for all $A\in\cal A$, is the GNS representation induced by $\omega$.
This simplifying assumption is similar to the assumption made in \cite{MP17} that
the probability distributions of the statistical manifold are absolutely
continuous w.r.t.~a given probability distribution.

\subsection*{Structure of the paper}

The next Section introduces the deformed logarithmic and exponential functions $\log_\phi$,
respectively $\exp_\phi$. In Section \ref {sect:states} the construction found in
\cite{MP17} is repeated with modifications to make it work in a non-commutative
context. The properties of the normalizing function are studied.
A class of states and their escorts is introduced.
A final Section gives a short discussion of the problems due to non-commutativity.

\section{The deformed logarithmic and exponential functions}
\label{sect:deffun}

\subsection{Definition and basic properties}

The specific deformed logarithmic and exponential functions introduced below
have been first used by Newton \cite {NNJ12}, without considering them as deformed functions.
The approach was then picked up by Montrucchio and Pistone \cite{MP17}.

Fix the function $\phi(u)=u/(1+u)$.  It is strictly positive and increasing on the interval $(0,+\infty)$. 
It determines a deformed logarithm \cite{NJ04} by
\be
\log_\phi(v)&=&\int_1^v\,\frac{1}{\phi(u)}\upd u\cr
&=&v-1+\log v.
\nnee
It is a concave function strictly increasing on $(0,+\infty)$.
The inverse function is denoted $\exp_\phi$. It is defined on all of the real line.
It is convex strictly increasing.
Special values are $\log_\phi(1)=0$ and $\exp_\phi(0)=1$.

Useful properties are
\be
\log_\phi(uv)&=&\log_\phi(u)+\log_\phi(v)+(u-1)(v-1),\\
\exp_\phi(u)&=&1+u-\log\exp_\phi(u),
\label {deform:id}\\
\frac{\upd\,}{\upd u}\exp_\phi(u)&=&\phi\left(\exp_\phi(u)\right)
=\frac{\exp_\phi(u)}{1+\exp_\phi(u)},\\
|\exp_\phi(u)-\exp_\phi(v)|&\le& |u-v|.
\label{form:ineq1}
\ee
The inequality follows from
\be
\exp_\phi(u)-\exp_\phi(v)\le u-v
\quad \mbox { if } u\ge v.
\label{form:ineq2}
\ee
To prove this use (\ref {deform:id}) and the fact that $\exp_\phi$  is an increasing function.
Inequality (\ref {form:ineq2}) also implies that 
\be
\exp_\phi(u)
&\le& 1+u\quad\mbox{if}\quad u\ge 0,\cr
&\ge&1+u\quad\mbox{if}\quad u\le 0.
\nnee
From $\log(1+u)\le u$ follows
\be
\log_\phi(1+u)=u+\log(1+u)\le 2u.
\nnee
This implies
\be
\exp_\phi(u)&\ge&1+\frac{u}{2}\quad\mbox{for all}\quad u.
\label{form:ineq4}
\ee
For $u<<0$ a better estimate follows from
\be
-u+\log\exp_\phi(u)&\ge&0
\qquad\mbox{ if }u\le 0.
\nnee
It implies
\be
e^{1+u}\ge\exp_\phi(u)\ge e^{u}
\qquad\mbox{ if }u<0.
\nnee

\subsection{Further properties}

The following results are needed later on.

\begin{proposition}
\label{form:prop:series}
 For all $u$ is $\exp_\phi(u)\le 1+\frac 12u+\frac 1{12} u^2$.
\end{proposition}

\beginproof
The function
\be
f(u)&=&1+\frac 12u+\frac 1{12} u^2-\exp_\phi(u).
\nnee
is convex. It attains its minimum at $u=0$ where its value equals 0.
\endproof

\begin{proposition}
\label{deform:prop:logsquareineq}
One has
\begin{description}
 \item 1) For all $u>0$ is
 \be
 \left[\log(1+u)\right]^2\le \left[\log 2\right]^2+(u-1)^2+\left[u-1+\log(1+u)\right]\log u.
 \label{escort:lemmastat1}
 \ee
 
 \item 2) There exists a constant $C$ such that
 \be
 \left[\log\,\phi(u)\right]^2\le C+(\log_\phi(u))^2\quad\mbox{ for all }u>0.
 \label{escort:lemmastat2}
 \ee
 
\end{description}
\end{proposition}
\beginproof

\paragraph{1)}
The derivative of r.h.s.~-~l.h.s.~equals
\be
& &
2(u-1)+\frac{1}{u}\left[[u-1+\log(1+u)\right]+\left[1+\frac{1}{1+u}\right]\log u
-\frac{2}{1+u}\log(1+u).\cr
& &
\nnee
It is a strictly increasing function which vanishes at $u=1$.
Hence, r.h.s.~-~l.h.s.~of (\ref {escort:lemmastat1}) is minimal at $u=1$.
One finally verifies that r.h.s.~=~l.h.s.~holds at $u=1$.

\paragraph{2)}
From 1) follows that
\be
\left[\log\phi(u)\right]^2
&=&
\left[\log(1+u)\right]^2+\left[\log u\right]^2-2\left[\log(1+u)\right]\left[\log u\right]\cr
&\le&
\left[\log 2\right]^2+(u-1)^2+\left[u-1+\log(1+u)\right]\log u\cr
& &
+\left[\log u\right]^2-2\left[\log(1+u)\right]\log u\cr
&=&
\left[\log 2\right]^2+(\log_\phi(u))^2+R(u)
\nnee
with
\be
R(u)&=&-\left[u-1+\log(1+u)\right]\log u.
\nnee
The function $R(u)$ is continuous and tends to $-\infty$ both when $u$ tends to 0 and when $u$ tends to $+\infty$.
It is positive on the interval $[u_0,1]$, with $u_0$ the solution of $\log(1+u_0)=1-u_0$.
In this interval it has a unique maximum, outside it is negative.
Hence the function $R(u)$ is bounded above. This implies the existence of a constant $C$ such that 
(\ref {escort:lemmastat2}) holds.

\endproof

Numerically, one finds that $C<0.52$ is feasible.

\begin{proposition}
\label{deform:prop:ftlambda}
Let the function $f_t$ be defined by
\be
f_t(\lambda,\mu)&=&\exp_\phi(t\log_\phi(\lambda)+\mu),
\nnee
for $t>0$, $\lambda>0$ and $\mu\in\Ro$.
\begin{itemize}
 \item [1)] $f_t(\lambda,\mu)$ is strictly increasing in each of the two arguments;
 \item [2)]
$0<t\le 1$ and $\mu+2(1-t)+t\log 2\ge 0$ implies
$\frac 12 t\lambda\le f_t(\lambda,\mu)$;
 \item [3)] $f_t(\lambda,\mu)\le t\lambda+\gamma$,
with  $\gamma=\exp\left(1+[\mu-t\log t]/[1-t]\right)$.
\end{itemize}
\end{proposition}

\beginproof
1)
The function $f_t$ is a composition of two strictly increasing functions.

2)
Let 
\be
g_t(\lambda,\mu)&=&\log_\phi(f_t(\lambda,\mu))-\log_\phi(t\lambda/2).
\nnee
At fixed value of $\mu\ge 0$ this function has a minimum when
$\lambda=2(1-t)/t$.
For this value of $\lambda$ is
\be
g_t(\lambda,\mu)&=&2(1-t)+\mu+t\log 2-t\log t-(1-t)\log (1-t)\cr
&\ge&2(1-t)+\mu+t\log 2.
\nnee
Hence, $g_t(\lambda,\mu)\ge 2(1-t)+\mu+t\log 2$ holds for all $\lambda$.
Because $\log_\phi$ is monotone increasing this implies the lower bound for 
the function $f_t$. 

3) Consider the function $h_t$ defined by
\be
h_t(\lambda,\mu)
&=&\log_\phi(t\lambda+\gamma)-\log_\phi(f_t(\lambda,\mu))\cr
&=&
\gamma-(1-t)+\log(t\lambda+\gamma)-t\log\lambda-\mu.
\nnee
It is minimal for $\lambda=\gamma/(1-t)$. For this value of $\lambda$ is
\be
h_t(\lambda,\mu)
&=&\gamma-t\log t -(1-t)\log(1-t)\ge 0.
\nnee
Because $\log_\phi$ is monotone increasing the upper bound for $f_t$ follows.

\endproof

\begin{proposition}
\begin{description}
 \item []
 \item [1)] For all $u$ is $[\exp_\phi(u)]^2\le 1+u\exp_\phi(u)$
 with equality if and only if $u=0$;
 \item [2)] The function $f$ defined by
 \be
 f(u)=\frac{1}{u}[\exp_\phi(u)-1]
 \nnee
 is increasing and satisfies $0<f(u)<1$;
 \item [3)] The derivative $f'$ of $f$ satisfies $f'<\frac 12$;
 \item [4)] For all $u,v$ is $ |f(u)-f(v)|\le \frac 12 |u-v|$.
\end{description}

\end{proposition}

\beginproof

\paragraph{1)}
Let $g(u)=1+u\exp_\phi(u)-[\exp_\phi(u)]^2$.
This function is strictly convex and has a minimum at $u=0$, with $g(0)=0$.

\paragraph{2)}
The derivative $f'$ of the function $f$ equals
\be
f'(u)=\frac{1}{u^2}\frac{1}{1+\exp_\phi(u)}g(u).
\nnee
Since $g(u)\ge 0$ for all $u$ the function $f$ is increasing.
A short calculation shows that $f(-\infty)=0$ and $f(+\infty)=1$.

\paragraph{3)}
The derivative $f'(u)$ goes through a maximum in two points.
They can be found by solving $f''(u)=0$.
For $u\not=0$ this equation is equivalent with
\be
u^2\exp_\phi(u)&=&2g(u)[1+\exp_\phi(u)]^2.
\nnee
At these points the expression for $f'(u)$ can be simplified to
\be
f'(u)=\frac{\exp_\phi(u)}{2[1+\exp_\phi(u)]^2}< \frac 12.
\nnee

\paragraph{4)}
This follows immediately from the previous items.

\endproof

\subsection{Note about operator monotonicity}

A function $f(u)$ is operator-monotone if $A\le B$ implies $f(A)\le f(B)$ for any pair
$A,B$ of self-adjoint matrices. 
The logarithm is operator-monotone \cite{PD08}. The deformed logarithm $\log_\phi(v)=v-1+\log v$
is the sum of operator-monotone functions. Hence it is also operator-monotone.
Finally, an oper\-ator-monotone function is also automatically operator-concave.
This implies that
\be
\log_\phi(\lambda A+(1-\lambda) B)\ge \lambda\log_\phi(A)+(1-\lambda)\log_\phi(B)
\ee
for any pair of positive bounded operators $A$ and $B$ and for any $\lambda\in[0,1]$.
See Section 11.6 of \cite{PD08}.
In the Appendix an example is given of a function which is increasing and concave
but not operator-monotone. This particular function gives useful results in the commutative
case \cite{MP17}, results which do not follow in the present non-commutative context.

\section{States and their escorts}
\label{sect:states}

\subsection{Construction}

Now follows the construction of a special class of self-adjoint operators $X$
all satisfying $X>0$ and $||X^{1/2}\Omega||=1$.

\begin{definition}
Given a vector state $\omega$ defined by the  normalized vector $\Omega$ in $\cal H$
and a self-adjoint operator $H$ with spectral decomposition
$H=\int\, \lambda\upd E_\lambda$ the $\omega$-expectation of $H$ is defined by
\be
\langle H\rangle_\omega&=&\int_{-\infty}^\infty\,\lambda\upd(E_\lambda\Omega,\Omega),
\nnee
provided that the integral converges absolutely.
This is the case when $\Omega$ belongs to the domain of $|H|^{1/2}$.
\end{definition}

Let $H=J|H|$ be the polar decomposition of $H$. Then one has
\be
\langle H\rangle_\omega&=&(J|H|^{1/2}\Omega,|H|^{1/2}\Omega).
\nnee
If $\Omega$ is in the domain of $H$ then one has
$\langle H\rangle_\omega=(H\Omega,\Omega)$.

\begin{proposition}
\label{prop:constr}
Let be given a self-adjoint operator $H$ on the Hilbert space ${\cal H}$.
Then one has
\begin{description}
 \item {1)} Any $\Psi$ in the domain of $|H|^{1/2}$
 also belongs to the domain of $[\exp_\phi(H-\beta)]^{1/2}$
for any real number $\beta$;
 \item {2)} The map $\beta\rightarrow\exp_\phi(H-\beta)$ is strongly continuous.
 \item {3)} The map $\beta\rightarrow ||[\exp_\phi(H-\beta)]^{1/2}\,\Psi||$ is strictly decreasing
 for any $\Psi\not=0$ in the domain of $H$.
\end{description}
If $\Omega$ is normalized, and in addition it belongs to the domain of $|H|^{1/2}$
and the expectation $\langle H\rangle_\omega$ vanishes,
then one has
\begin{description}
 \item {4)} The map $\beta\rightarrow ||[\exp_\phi(H-\beta)]^{1/2}\,\Omega||$
 is one-to one from $(-\infty,+\infty)$  to $(0,+\infty)$.
 \item {5)} There exists a unique non-negative number, denoted $\alpha(H)$, for which
$||[\exp_\phi(H-\alpha(H))]^{1/2}\,\Omega||=1$.
\item {6)} For any real number $c$ is $\alpha(H+c)=\alpha(H)+c$.

\end{description}

\end{proposition}

\beginproof\noindent

\paragraph{1)}
From (\ref {form:ineq1}) follows that
\be
\exp_\phi(\lambda-\beta)&\le&\exp_\phi(-\beta)+|\lambda|.
\nnee
This implies 
\be
\exp_\phi(H-\beta)&\le&\exp_\phi(-\beta)+|H|
\nnee
and
\be
||[\exp_\phi(H-\beta)]^{1/2}\,\Psi||^2
&\le&\exp_\phi(-\beta)||\Psi||^2+||\,|H|^{1/2}\Psi||^2\cr
&<& +\infty.
\nnee
This shows that $\Psi$ belongs to the domain of $[\exp_\phi(H-\beta)]^{1/2}$.

\paragraph {2)}
Let 
\be
H=\int\, \lambda\upd E_\lambda
\nnee
denote the spectral decomposition of the operator $H$.
From (\ref {form:ineq1}) follows
\be
\exp_\phi(H-\beta_1)-\exp_\phi(H-\beta_2)
&=&
\int\,\left[\exp_\phi(\lambda-\beta_1)-\exp_\phi(\lambda-\beta_2)\right]\upd E_\lambda\cr
&\le&\int\,|\beta_1-\beta_2|\upd E_\lambda\cr
&=&|\beta_1-\beta_2|.
\nnee
This implies strong continuity of the map
\be
\beta\rightarrow \exp_\phi(H-\beta).
\nnee

\paragraph{3)}
Assume $\gamma>\beta$.
From the convexity of $\exp_\phi$ follows that 
\be
\exp_\phi(\lambda-\beta)
&\ge&
\exp_\phi(\lambda-\gamma)+(\gamma-\beta)\phi(\exp_\phi(\lambda-\gamma)).
\nnee
This implies
\be
& &||[\exp_\phi(H-\beta)]^{1/2}\,\Psi||^2\cr
&=&
\int\exp_\phi(\lambda-\beta)\upd (E_\lambda\Psi,\Psi)\cr
&\ge&
\int\left[\exp_\phi(\lambda-\gamma)+(\gamma-\beta)\phi(\exp_\phi(\lambda-\gamma))\right]\upd (E_\lambda\Psi,\Psi)\cr
&=&
||[\exp_\phi(H-\gamma)]^{1/2}\,\Psi||^2
+(\gamma-\beta)||[\phi(\exp_\phi(H-\gamma))]^{1/2}\,\Psi||^2.\cr
& &
\nnee
Because $\phi$ is strictly positive zero cannot be an eigenvalue of $\phi(\exp_\phi(H-\gamma))$.
Hence, $|[\phi(\exp_\phi(H-\gamma))]^{1/2}\,\Psi||\not=0$ and $\gamma>\beta$ implies the strict inequality
\be
||[\exp_\phi(H-\beta)]^{1/2}\Psi||^2&>&||[\exp_\phi(H-\gamma)]^{1/2}\Psi||^2.
\nnee

\paragraph{4)}
Introduce the notation $X_\beta(H)\equiv \exp_\phi(H-\beta)$.
Because $\log_\phi$ is concave and $\langle H\rangle_\omega=0$ one has
\be
-\beta
&=&\int_{-\infty}^{+\infty}(\lambda-\beta)\,\upd(E_\lambda\Omega,\Omega)\cr
&=&\langle\log_\phi(X_\beta(H)\rangle_\omega\cr
&\le&\log_\phi(||(X_\beta(H)]^{1/2}\,\Omega||^2.
\nnee
Hence, if $\beta$ tends to $-\infty$ then $\log_\phi(||(X_\beta(H)]^{1/2}\,\Omega||$ tends to $+\infty$.
This implies that $||(X_\beta(H)]^{1/2}\,\Omega||$ tends to $+\infty$.

On the other hand, if $\beta$ tends to $+\infty$ then $||[X_\beta(H)]^{1/2}\,\Omega||$ tends to zero.
This follows from the following argument.

Fix $\epsilon>0$. Because $\Omega$ is in the domain of $|H|^{1/2}$
there exists $\lambda_\epsilon>>0$ such that
\be
\int_{\lambda_\epsilon}^{+\infty}\lambda\upd (E_\lambda\Omega,\Omega)<\epsilon.
\nnee
Without restriction take $\lambda_\epsilon\ge 1$.
Next choose $\beta_\epsilon$ large enough so that
\be
\exp_\phi(\lambda_\epsilon-\beta_\epsilon)<\epsilon.
\nnee
Then for all $\beta>\beta_\epsilon$ is
\be
||[X_\beta(H)]^{1/2}\,\Omega||^2
&=&
\int_{-\infty}^{\lambda_\epsilon}\exp_\phi(\lambda-\beta)\upd (E_\lambda\Omega,\Omega)\cr
& &
+\int_{\lambda_\epsilon}^{+\infty}\exp_\phi(\lambda-\beta)\upd (E_\lambda\Omega,\Omega)\cr
&\le&
\epsilon\int_{-\infty}^{\lambda_\epsilon}\upd (E_\lambda\Omega,\Omega)\cr
& &
+\int_{\lambda_\epsilon}^{+\infty}\left[\exp_\phi(-\beta)+\lambda\right]\upd (E_\lambda\Omega,\Omega)\cr
&\le&
\epsilon+2\int_{\lambda_\epsilon}^{+\infty}\lambda\upd (E_\lambda\Omega,\Omega)\cr
&\le&3\epsilon.
\nnee
This finishes the proof that $||[X_\beta(H)]^{1/2}\Omega||$ tends to 0 as $\beta$ tends to $+\infty$.

Because $\beta\rightarrow ||[X_\beta(H)]^{1/2}\,\Omega||$ is strictly decreasing and continuous
one concludes that the map is one-to one.

\paragraph{5)}
The existence of a unique real number $\alpha(H)$ is an immediate consequence of item 4).
Convexity of the deformed exponential implies
\be
||[\exp_\phi(H)]^{1/2}\Omega||^2
&=&\int\,\exp_\phi(\lambda)\upd (E_\lambda\Omega,\Omega)\cr
&\ge&\exp_\phi\left(\int\upd (E_\lambda\Omega,\Omega)\,\lambda\right)\cr
&=&\exp_\phi\left(\langle H\rangle_\omega\right)\cr
&=&\exp_\phi(0)\cr
&=&1.
\nnee
Because the map $\beta\rightarrow ||[\exp_\phi(H-\beta)]^{1/2}\,\Omega||$ is strictly decreasing
and has a value $\ge 1$ at $\beta=0$ one concludes that it takes the value 1 at $\alpha(H)\ge0$.

\paragraph{6)}
This follows immediately from
\be
||\left[\exp_\phi([H+c]-[\alpha(H)+c])\right]^{1/2}\Omega||
&=&
||\left[\exp_\phi(H-\alpha(H))\right]^{1/2}\Omega||\cr
&=&1.
\nnee

\endproof

In the commutative case the function
$\alpha(H)$ is shown to be convex. See Proposition 1 of \cite{MP17}.
No such result is expected here because
unbounded self-adjoint operators do not form an affine space.
In addition, the convexity proof is based on the convexity of the
deformed exponential function. However, a convex increasing function
cannot be operator-monotone because operator-monotone functions
are automatically concave. It is therefore not immediately clear how
to prove convexity properties for these functions of operators.

\subsection{Properties of the normalization function}
\label{sect:propnorm}

Some properties of the normalization function $\alpha(H)$
are gathered in the following Proposition. Differentiability is considered below
in  Section \ref{sect:tangent} on tangent vectors.

\begin{proposition}
\label{constr:prop:deriv}
Let $H$ be a self-adjoint operator such that $\Omega$ belongs to the domain of $|H|^{1/2}$
and $\langle H\rangle_\omega=0$.
Let $\alpha(H)$ be the function defined in Proposition \ref{prop:constr}.
One has
\begin{description}
\item [1)] The function \, $t\in\Ro\mapsto \alpha(tH)$ \, is convex;
\item [2)] $\alpha(H)\ge 0$;
\end{description}
If in addition $\Omega$ is in the domain of $H$ then also the following holds.
\begin{description}
\item [3)] $||H\Omega||<1$ \, implies \, $\alpha(H)<||H\Omega||^2$;
In particular, the derivative of
the function $t\in\Ro\mapsto \alpha(tH)$ exists and vanishes at $t=0$.
\item [4)] For any $\beta$ is
\be
||\,[\exp_\phi(H-\beta)]^{1/2}\Omega||^2
&\le& 1-\frac 12\beta+\frac 1{12}\left(||H\Omega||^2+\beta^2\right);
\nnee
\item [5)] $\alpha(H)\ge 6-\sqrt{36-||H\Omega||^2}\ge\frac 1{12}||H\Omega||^2$;
\item [6)] If \, $||H\Omega||\le 3$ \, then \, $\alpha(H)\le ||H\Omega||$.
\end{description}

\end{proposition}

\beginproof

\paragraph{1)}
Let
\be
H&=&\int\lambda\upd E_\lambda
\nnee
be the spectral decomposition of $H$.
Let
\be
\beta=\mu\alpha(t_1 H)+(1-\mu)\alpha(t_2 H).
\nnee
One has
\be
& &
||[\exp_\phi([\mu t_1+(1-\mu) t_2]H-\beta)]^{1/2}\Omega||^2\cr
&=&
\int\upd (E_\lambda\Omega,\Omega)\,\exp_\phi ((\mu t_1+(1-\mu) t_2)\lambda-\beta)\cr
&\le&
\mu\int\upd (E_\lambda\Omega,\Omega)\,\exp_\phi (t_1\lambda-\alpha(t_1 H))\cr
& &+(1-\mu)\int\upd (E_\lambda\Omega,\Omega)\,\exp_\phi ( t_2\lambda-\alpha(t_2 H))\cr
&=&
\mu ||[\exp_\phi(t_1H-\alpha(t_1 H))]^{1/2}\Omega||^2
+(1-\mu)||[\exp_\phi(t_2H-\alpha(t_2 H))]^{1/2}\Omega||^2\cr
&=&1.
\nnee
This implies that $\beta\ge \alpha((\mu t_1+(1-\mu) t_2)H)$.
One concludes that $\mu\alpha(t_1 H)+(1-\mu)\alpha(t_2 H)\ge \alpha((\mu t_1+(1-\mu) t_2)H)$,
i.e.\, $t\mapsto \alpha(tH)$ \, is convex.

\paragraph{2)}
With the help of the inequality (\ref {form:ineq4}) it follows that
\be
||[\exp_\phi(H)]^{1/2}\Omega||^2
&=&
\int\upd (E_\lambda\Omega,\Omega)\,\exp_\phi (\lambda)\cr
&\ge&
\int\upd (E_\lambda\Omega,\Omega)\,\left(1+\frac 12\lambda\right)\cr
&=&
1+\frac 12\langle H\rangle_\omega\cr
&=&1.
\nnee
This implies $0\le\alpha(H)$.

\paragraph{3)}

Because $\alpha(tH)\ge 0$ for all $t$, $\alpha(0)=0$ and $t\mapsto \alpha(tH)$ is convex
one concludes that $\alpha(tH)$ is continuous with a minimum at $t=0$. 
Now calculate, using (\ref {deform:id}) and $\langle H\rangle_\omega=0$,
\be
& &
||[\exp_\phi(tH-\gamma t^2)]^{1/2}\Omega||^2\cr
&=&
\int\upd (E_\lambda\Omega,\Omega)\,\exp_\phi (t\lambda-\gamma t^2)\cr
&=&
-\frac 12\gamma t^2
+\int\upd (E_\lambda\Omega,\Omega)\,\left[
\exp_\phi (t\lambda-\gamma t^2)-\frac 12(t\lambda-\gamma t^2)\right]\cr
&=&
1-\frac 12\gamma t^2\cr
& &
+\int\upd (E_\lambda\Omega,\Omega)\,\left[
\frac 12(t\lambda-\gamma t^2)-\log \exp_\phi (t\lambda-\gamma t^2)\right].
\nnee
Split the integral in two using the constant $\mu=\gamma t-1/t$. This gives
\be
& &
||[\exp_\phi(tH-\gamma t^2)]^{1/2}\Omega||^2\cr
&=&
1-\frac 12\gamma t^2\cr
& &
+\int_{-\infty}^\mu\upd (E_\lambda\Omega,\Omega)\,\left[
\frac 12(t\lambda-\gamma t^2)-\log \exp_\phi (t\lambda-\gamma t^2)\right]\cr
& &
+\int_\mu^{+\infty}\upd (E_\lambda\Omega,\Omega)\,\left[
\frac 12(t\lambda-\gamma t^2)-\log \exp_\phi (t\lambda-\gamma t^2)\right].
\nnee
In the last term use (\ref {form:ineq4}) and $\log(1+u)\ge u-u^2$ when $u>-1/2$ to obtain
\be
& &
||[\exp_\phi(tH-\gamma t^2)]^{1/2}\Omega||^2\cr
&\le&
1-\frac 12\gamma t^2\cr
& &
+\int_{-\infty}^\mu\upd (E_\lambda\Omega,\Omega)\,\left[
\frac 12(t\lambda-\gamma t^2)-\log \exp_\phi (t\lambda-\gamma t^2)\right]\cr
& &
+\int_\mu^{+\infty}\upd (E_\lambda\Omega,\Omega)\,\left[
\frac 12(t\lambda-\gamma t^2)-\log(1+\frac 12(t\lambda-\gamma t^2))\right]\cr
&\le&
1-\frac 12\gamma t^2\cr
& &
+\int_{-\infty}^\mu\upd (E_\lambda\Omega,\Omega)\,\left[
\frac 12(t\lambda-\gamma t^2)-\log \exp_\phi (t\lambda-\gamma t^2)\right]\cr
& &
+\int_\mu^{+\infty}\upd (E_\lambda\Omega,\Omega)\,\left[
\frac 14(t\lambda-\gamma t^2)^2\right]\cr
&=&
1-\frac 12\gamma t^2+\frac 14 t^2||(H-\gamma t)\Omega||^2\cr
& &
+\int_{-\infty}^\mu\upd (E_\lambda\Omega,\Omega)\,\left[
\frac 12(t\lambda-\gamma t^2)-\log \exp_\phi (t\lambda-\gamma t^2)
-\frac 14(t\lambda-\gamma t^2)^2
\right].
\nnee
Note that $\log\exp_\phi(u)\ge u$ holds for all $u\le 0$. 
Hence,
\be
& &
||[\exp_\phi(tH-\gamma t^2)]^{1/2}\Omega||^2\cr
&\le&
1-\frac 12\gamma t^2+\frac 14 t^2||(H-\gamma t)\Omega||^2
+\int_{-\infty}^\mu\upd (E_\lambda\Omega,\Omega)\,f(\lambda),
\nnee
with the function $f(\lambda)$ defined by
\be
f(\lambda)=-\frac 12(t\lambda-\gamma t^2)-\frac 14(t\lambda-\gamma t^2)^2.
\nnee
This function is increasing on the interval $(-\infty,\mu]$.
Using $f(\lambda)\le f(\mu)=-1$ one obtains
\be
& &
||[\exp_\phi(tH-\gamma t^2)]^{1/2}\Omega||^2\cr
&\le&
1-\frac 12\gamma t^2+\frac 14 t^2||(H-\gamma t)\Omega||^2
+\int_{-\infty}^\mu\upd (E_\lambda\Omega,\Omega)\,f(\mu)\cr
&=&
1-\frac 12\gamma t^2+\frac 14 t^2||(H-\gamma t)\Omega||^2
-(E_\mu\Omega,\Omega).
\nnee

Take now $\gamma=||H\Omega||^2$.
Then one obtains
\be
||[\exp_\phi(tH-\gamma t^2)]^{1/2}\Omega||^2
\le
1-\frac 14\gamma t^2(1-\gamma t^2).
\nnee
This shows that $||[\exp_\phi(tH-\gamma t^2)]^{1/2}\Omega||<1$ for $|t|<1/\sqrt\gamma$.
The latter implies that $\alpha(tH)<\gamma t^2$ for $|t|<1/\sqrt\gamma$.
If now $||H\Omega||<1$ holds then one can take $t=1$ to obtain $\alpha(H)<||H\Omega||^2$.

\paragraph{4)}
Apply the inequality of Proposition \ref {form:prop:series} to obtain
\be
||\,[\exp_\phi(H-\beta)]^{1/2}\Omega||^2
&=&
\int_{-\infty}^{+\infty}\exp_\phi(\lambda-\beta)\,\upd(E_\lambda\Omega,\Omega)\cr
&\le&
\int_{-\infty}^{+\infty}\left[ 1+\frac 12(\lambda-\beta)+\frac 1{12} (\lambda-\beta)^2\right]
\,\upd(E_\lambda\Omega,\Omega)\cr
&=&
1+\frac 12(\langle H\rangle_\omega-\beta)+
\frac 1{12}\left(||H\Omega||^2-2\beta \langle H\rangle_\omega+\beta^2\right)\cr
&=&
1-\frac 12\beta+\frac 1{12}\left(||H\Omega||^2+\beta^2\right).
\nnee

\paragraph{5)}

One has $\beta\le\alpha(H)$ if and only if $||\,[\exp_\phi(H-\beta)]^{1/2}\Omega||^2\ge 1$.
The inequality proved above shows that a sufficient condition is that 
$\beta^2-6\beta+||H\Omega||^2\ge 0$. The zeroes of this quadratic equation are
$\beta=6\pm\sqrt{36-||H\Omega||^2}$. Hence, $\beta\le 6-\sqrt{36-||H\Omega||^2}$
suffices to obtain $\beta\le\alpha(H)$.

\paragraph{6)}
Calculate, using Proposition \ref {form:prop:series},
\be
||[\exp_\phi(H-3)]^{1/2}\Omega||^2
&=&
\int_{-\infty}^{+\infty}\exp_\phi(\lambda-3)\,\upd(E_\lambda\Omega,\Omega)\cr
&\le&
\int_{-\infty}^{+\infty}
\left[1+\frac 12(\lambda-3)+\frac 1{12}(\lambda-3)^2
\right]
\,\upd(E_\lambda\Omega,\Omega)\cr
&=&
\frac 14+\frac 1{12}||H\Omega||^2\cr
&\le&1.
\nnee
This implies that $\alpha(H)\le 3$. Because $\alpha(0)=0$ and the function $\alpha$
is convex one concludes that $\alpha(H)\le ||H\Omega||$ as long as $||H\Omega||\le 3$. 

\endproof

\subsection{States}

\begin{definition}
Given a positive operator $X$, which satisfies $||X^{1/2}\Omega||=1$,
let $\omega_X$ denote the vector state of $\cal A$ defined by
\be
\omega_X(A)&=&
(A X^{1/2}\Omega, X^{1/2}\Omega)\quad\mbox{ for all }A\in {\cal A}.
\label {def:omegaXdef}
\ee
\end{definition}
Note that $\omega_X=\omega_Y$ does not necessarily imply that $X=Y$.

\begin{proposition}
\label{def:prop:core}
Let be given a self-adjoint operator $H$ which is affiliated with 
the commutant ${\cal A}'$.
Assume $\Omega$ is in the domain of $|H|^{1/2}$ and $\langle H\rangle_\omega=0$.
Let $X=\exp_\phi(H-\alpha(H))$, with $\alpha(H)$ the function defined by Proposition \ref {prop:constr}.
The following holds:
\begin{itemize}
 \item [1)] The domain of $|H|^{1/2}$ is a subspace of the domain of $X^{1/2}$;
 \item [2)] Any $\psi$ in $\dom |H|^{1/2}$
 satisfies
 \be
 ||X^{1/2}\psi||^2&\le& ||\psi||^2+||\,|H|^{1/2}\psi||^2;
 \nnee
 \item [3)] $X^{1/2}\Omega$ is
 separating for $\cal A$;
 \item [4)] If ${\cal A}\Omega$ is a core of $|H|^{1/2}$
 then ${\cal A}\Omega$ is a core of $X^{1/2}$ as well;
 \item [5)] If ${\cal A}\Omega$ is a core of $X^{1/2}$ then $X^{1/2}\Omega$ is cyclic for $\cal A$.
\end{itemize}

\end{proposition}

\beginproof

\paragraph{1)}
This follow immediately from Proposition \ref {prop:constr}.

\paragraph{2)}
Let
\be
H&=&\int_{-\infty}^{+\infty}\lambda\,\upd E_\lambda
\nnee
denote the spectral decomposition of $H$.
Take $\psi$ in the domain of $|H|^{1/2}$.
From the spectral theorem it follows that
$\psi$ belongs to the domain of $X^{1/2}$ if and only if
\be
||X^{1/2}\psi||^2&=&
\int_{-\infty}^{+\infty}\exp_\phi(\lambda-\alpha(H))\,\upd (E_\lambda\psi,\psi)
\nnee
remains finite.
From $\alpha(H)\ge 0$ and the inequality $\exp_\phi(u)\le 1+\max\{0,u\}$
one obtains
\be
||X^{1/2}\psi||^2&\le&
\int_{-\infty}^{+\infty}[1+\max\{0,\lambda-\alpha(H)\})]\,\upd (E_\lambda\psi,\psi)\cr
&=&
\int_{-\infty}^{\alpha(H)}\,\upd (E_\lambda\psi,\psi)\cr
& &
+\int_{\alpha(H)}^{+\infty}[1+\lambda-\alpha(H))]\,\upd (E_\lambda\psi,\psi)\cr
&\le&
||\psi||^2+\int_{\alpha(H)}^{+\infty}\lambda\,\upd (E_\lambda\psi,\psi)\cr
&\le&
||\psi||^2+||\,|H|^{1/2}\psi||^2\cr
&<&+\infty.
\nnee
This proves 2).

\paragraph{3)}
Take $A\in{\cal A}$ and assume $AX^{1/2}\Omega=0$.
Because $A$ and $X^{1/2}$ commute this implies $X^{1/2}A\Omega=0$.
The operator $X^{1/2}$ is strictly positive and therefore invertible.
Hence, it follows that $A\Omega=0$. However, $\Omega$ is separating for $\cal A$.
Therefore, one concludes that $A=0$. 
This shows that $X^{1/2}\Omega$ is separating for $\cal A$.

\paragraph{4)}
Take $\psi$ in the domain of $X^{1/2}$.
Let $\psi_n=(E_n-E_{-n})\psi$. These vectors
belong to the domain of $|H|^{1/2}$ because $|H|^{1/2}\psi_n=[|H|^{1/2}(E_n-E_{-n})]\psi_n$
and $|H|^{1/2}(E_n-E_{-n})$ is a bounded operator.
Because ${\cal A}\Omega$ is a core of $|H|^{1/2}$ there exist $A_{n,m}$ in $\cal A$
such that
\be
\psi_n=\lim_mA_{n,m}\Omega
\quad\mbox{ and }\quad
|H|^{1/2}\psi_n=\lim_m |H|^{1/2}A_{n,m}\Omega.
\nnee 
Use $\alpha(H)\ge 0$ and the inequality $\exp_\phi(u)\le 1+\max\{0,u\}$
to obtain
\be
& &
||X^{1/2}[A_{n,m}-A_{N,p}]\Omega||^2\cr
&\le&
||[A_{n,m}-A_{N,p}]\Omega||^2\cr
& &
+\int_{\alpha(H)}^{+\infty}[\lambda-\alpha(H)]\upd(E_\lambda [A_{n,m}-A_{N,p}]\Omega,[A_{n,m}-A_{N,p}]\Omega)\cr
&\le&
||[A_{n,m}-A_{N,p}]\Omega||^2+||\,|H|^{1/2}[A_{n,m}-A_{N,p}]\Omega||^2.
\nnee
This shows that the vectors $X^{1/2}A_{n,m}\Omega$ form a Cauchy sequence.
Because $X^{1/2}$ is a closed operator it necessarily converges to $X^{1/2}\psi_n$.
It is then easy to show that
\be
\lim_n\lim_m X^{1/2}A_{n,m}\Omega=\lim_n X^{1/2}\psi_n=X^{1/2}\psi.
\nnee
One concludes that ${\cal A}\Omega$ is a core of $X^{1/2}$.

\paragraph{5)}
Assume $\psi$ is orthogonal to ${\cal A}X^{1/2}\Omega$.
Take $\phi$ in the domain of $X^{1/2}$.
By assumption is ${\cal A}\Omega$ a core of $X^{1/2}$. Hence there exist $A_n$ in $\cal A$
such that $A_n\Omega$ converge to $\phi$ and $X^{1/2}A_n\Omega$ converge to $X^{1/2}\phi$.
From $(X^{1/2} A_n \Omega, \psi)=(A_n X^{1/2}\Omega, \psi)=0$ then follows that 
$(X^{1/2}\phi, \psi)=0$. Because $X$ is invertible the range of $X$ is dense in $\cal H$.
Therefore one concludes that $\psi=0$.

\endproof

\begin{proposition}
Let $X$ and $Y$ satisfy the conditions of Proposition \ref {def:prop:core}.
Assume ${\cal A}\Omega$ is a core of $X^{1/2}$ and of $Y^{1/2}$.
Then $\omega_X=\omega_Y$ implies that $X=Y$.
\end{proposition}

\beginproof
From $\omega_X=\omega_Y$ follows that there exists an isometry $U$ in ${\cal A}'$
such that $Y^{1/2}\Omega=UX^{1/2}\Omega$. Because $X$ and $Y$ are affiliated with the
commutant ${\cal A}'$ this implies that $Y^{1/2}=UX^{1/2}$ holds
on the space ${\cal A}\Omega$. The latter is by assumption a core of $X^{1/2}$.
One concludes that $Y^{1/2}=UX^{1/2}$ on the domain of $X^{1/2}$. Because of the uniqueness
of the polar decomposition this implies that $U=1$ and $Y=X$.

\endproof

\subsection{Escort states}
\begin{proposition}
\label{escort:def:prop1}
Let $X$ be a strictly positive operator satisfying $||X^{1/2}\Omega||=1$. One has
\begin{description}
 \item 1) $0<\phi(X)<1$; in particular, $\phi(X)$ is a bounded strictly positive operator with norm
 $||\phi(X)||\le 1$.
 \item 2) $0<(\phi(X)\Omega,\Omega)\le 1/2$;
 \item 3) $\phi(X)$ belongs to the commutant ${\cal A}'$.
\end{description}
\end{proposition}

\beginproof
\paragraph{1)}
This follows because $X>0$ and $0<\phi(u)<1$ for all $u>0$.

\paragraph{2)}
Use the spectral decomposition of $X$
\be
X&=&\int_0^{+\infty}\lambda\upd E_\lambda,
\nnee
together with the concavity of the function $\phi$ and
the normalization\hfil \\
$||X^{1/2}\Omega||^2=1$ to write
\be
(\phi(X)\Omega,\Omega)
&=&
\int_0^{+\infty}\phi(\lambda)\upd (E_\lambda\Omega,\Omega)\cr
&\le&
\phi\left(\int_0^{+\infty}\lambda\upd (E_\lambda\Omega,\Omega)\right)\cr
&=&
\phi\left(||X^{1/2}\Omega||^2\right)\cr
&=&\phi(1)\cr
&=&\frac 12.
\nnee
Finally $(\phi(X)\Omega,\Omega)=0$ implies $||X^{1/2}\Omega||=0$, which 
contradicts $||X^{1/2}\Omega||=1$. One concludes that 
$0<(\phi(X)\Omega,\Omega)\le 1/2$.

\paragraph{3)}
That $\phi(X)$ belongs to ${\cal A}'$ follows because it is a bounded function of
a self-adjoint operator affiliated with ${\cal A}'$.

\endproof

\begin{definition}
Given a self-adjoint positive operator $X$ affiliated with the commutant ${\cal A}'$
and satisfying $||X^{1/2}\Omega||=1$
introduce the state $\tilde\omega_X$
of $\cal A$ defined by
\be
\tilde\omega_X(A)&=&\frac{(A\Omega,\phi(X)\Omega)}{(\Omega,\phi(X)\Omega)},
\qquad A\in {\cal A}.
\nnee
The state $\tilde\omega_X$ is an {\em escort} of $\omega_X$.
\end{definition}

\subsection{Tangent vectors}
\label{sect:tangent}

Escort states appear in a natural manner when studying vectors tangent
to geodesics.

\begin{proposition}
\label{constr:prop:diff}
Let $H$ be a self-adjoint operator such that $\Omega$ belongs to the domain of $|H|^{1/2}$
and $\langle H\rangle_\omega=0$.
Let $\alpha(H)$ be the function defined in Proposition \ref{prop:constr}.
The function $t\in\Ro\mapsto \alpha(tH)$ is differentiable.
The derivative satisfies
\be
\frac{\upd\,}{\upd t}\alpha(tH)
&=&\frac{(J|H|^{1/2}\Omega,\phi(X_t)|H|^{1/2}\Omega)}{(\Omega,\phi(X_t)\Omega)},
\nnee
with $X_t$ given by $X_t=\exp_\phi(tH-\alpha(tH))$.
\end{proposition}

\beginproof

From the identity (\ref {deform:id}) it follows that
\be
1&=&||X_t^{1/2}\Omega||^2\cr
&=&\int \exp_\phi(t\lambda-\alpha(tH))\,\upd(E_\lambda\Omega,\Omega)\cr
&=&\int \left[1+t\lambda-\alpha(tH)
-\log\exp_\phi(t\lambda-\alpha(tH))\right]\,\upd(E_\lambda\Omega,\Omega)\cr
&=&
1-\alpha(tH)
-\int \log\exp_\phi(t\lambda-\alpha(tH))\,\upd(E_\lambda\Omega,\Omega).
\nnee
This can be written as
\be
\alpha(tH)
&=&
-\int f(t\lambda-\alpha(tH))\,\upd(E_\lambda\Omega,\Omega).
\nnee
with $f(u)=\log\exp_\phi(u)$. The function $f(u)$ is concave. This implies
\be
f(u)-f(v)&\le&(u-v)f'(v),\qquad u,v\in\Ro.
\nnee
Hence one obtains
\be
& &\alpha(t'H)-\alpha(tH)\cr
&=&
-\int [f(t'\lambda-\alpha(t'H))-f(t\lambda-\alpha(tH))]\,\upd(E_\lambda\Omega,\Omega)\cr
&\ge&
-\int [(t'-t)\lambda-\alpha(t'H)+\alpha(tH)]f'(t\lambda-\alpha(tH))\,\upd(E_\lambda\Omega,\Omega).
\nnee
Use $f'(u)=1/(1+\exp_\phi(u))$ to show that for any real $t,t'$ one has the inequality
\be
\alpha(t'H)-\alpha(tH)
&\ge&
(t'-t) f(t)
\nnee
with
\be
f(t)&=&\frac{\int \lambda\phi(\exp_\phi(t\lambda-\alpha(tH)))\,\upd(E_\lambda\Omega,\Omega)}
{\int \phi(\exp_\phi(t\lambda-\alpha(tH)))\,\upd(E_\lambda\Omega,\Omega)}.
\nnee
Swap $t$ and $t'$ to obtain
\be
(t'-t)f(t')\ge \alpha(t'H)-\alpha(tH) \ge (t'-t) f(t).
\nnee
Because the function $f(t)$
is continuous one concludes
that it is the derivative of $t\mapsto \alpha(tH)$.

\endproof

Let $\omega_t\equiv\omega_{X_t}$ with $X_t$ defined as in the above Proposition.
A short calculation shows that
\be
\frac{\upd\,}{\upd t}\omega_t
&=&f_t
\nnee
with for any $A\in{\cal A}$
\be
f_t(A)&=&
(JA|H|^{1/2}\Omega,\phi(X_t)|H|^{1/2}\Omega)-(A\Omega,\phi(X_t)\Omega)\frac{\upd\,}{\upd t}\alpha(tH).
\nnee
The linear functional $f_t$ belongs to the dual space ${\cal A}^*$
and is a vector tangent to the curve $t\mapsto\omega_t$.

%%%%%%%%%%%%%%%%%%%%%%%%%%%%%%%%%%%%%%
\section{Discussion}

Part of the work of Montrucchio and Pistone \cite{MP17} is transferred to a non-commutative
setting in a rather straightforward manner. The probability distributions are replaced by
vector states on a von Neumann algebra. Probability densities are replaced by positive
operators affiliated with the commutant of the von Neumann algebra.
The properties of the normalization function are studied
in Section \ref {sect:propnorm}.

The main obstacle in generalizing all of \cite{MP17} to a non-commutative context
is that certain monotone functions and convex functions, appearing in the proofs of \cite{MP17},
are not operator-monotone, respectively operator-convex. See the Appendix below.
In addition, 
technical difficulties arise because the sum of two self-adjoint operators is in general
not self-adjoint due to problems with the domain of definition.
These difficulties prevent a straightforward introduction of a geometric structure
on the manifold of faithful vector states.

%%%%%%%%%%%%%%%%%%%%%%%%%%%%%%%%%%%%%%%
\appendix
\section*{Appendix}

The following negative results give an indication of the kind
of problems that one encounters with functions of operators.
In this Appendix the function $\phi$ is defined by
\be
\phi(u)&=&\frac{u}{\lambda+u},
\nnee
where $\lambda$ is a fixed positive constant.
The deformed logarithm equals
\be
\log_\phi(u)&=&u-1+\lambda\log u.
\nnee

\paragraph{Proposition}
The function $f(u)=u-\exp_\phi(u)$ is {\em not} operator-monotone.
\vskip 0.1cm

\beginproof

Introduce the shorthands $x=\exp_\phi(u)$ and $y=\exp_\phi(v)$.
One has
\be
f'(u)
&=&1-\phi(\exp_\phi(u))\cr
&=&1-\frac{x}{\lambda+x}\cr
&=&\frac{\lambda}{\lambda+x}.
\nnee
A necessary condition (Theorem 11.17 of \cite{PD08}) for $f(u)$ to be operator-monot\-one
is that the following determinant is positive
\be
D
&=&
\left|\begin{array}{lr}
       f'(u) &1-\frac{x-y}{u-v}\\
       1-\frac{x-y}{u-v} &f'(v)
      \end{array}\right|\cr
&=&\frac{\lambda}{\lambda+x}\,\frac{\lambda}{\lambda+y}
-\left[1-\frac{x-y}{u-v}\right]^2.
\nnee
Consider the case $u>v$. This implies $x>y$.
Introduce $\epsilon>0$ defined by $x=(1+\epsilon)y$.
Then
\be
u
&=&\log_\phi(x)\cr
&=&x-1+\lambda\log x\cr
&=&y-1+\lambda\log y +\epsilon y+\lambda\log(1+\epsilon)\cr
&=&v+\epsilon y+\lambda\log(1+\epsilon).
\nnee
One obtains
\be
D
&=&
\frac{\lambda}{\lambda+y+\epsilon y}\,\frac{\lambda}{\lambda+y}
-\left[\frac{\lambda\log(1+\epsilon)}{\epsilon y+\lambda\log(1+\epsilon)}\right]^2.
\nnee
The condition that $D>0$ becomes
\be
(\lambda+y)(\lambda+y+\epsilon y)\left[\frac{\log(1+\epsilon)}{\epsilon y+\lambda\log(1+\epsilon)})\right]^2<1
\quad\mbox{for all }y>0,\epsilon>0.
\nnee
Let $\delta>1$ be given by
\be
\delta=\frac{\epsilon}{\log(1+\epsilon)}.
\nnee
Then this condition becomes
\be
(\lambda+y)(\lambda+y+\epsilon y)<[\lambda+\delta y]^2,
\nnee
or, equivalently, 
\be
[\delta^2-(1+\epsilon)]y>\lambda\left[\epsilon-2(\delta-1)\right].
\ee
This equation puts a condition on the choice of $y$.
Take for instance $\epsilon=\delta=e-1$. Then the condition reads
\be
y>\lambda\frac{3-e}{e^2-3e+1}.
\nnee
The r.h.s.~of this condition is positive. Hence there exist choices of $y>0$
which do not satisfy the condition.

\endproof

\paragraph{Corollary}
There exist hermitian matrices
$A$ and $B$ which violate the operator version of (\ref {form:ineq2}),
i.e.~for which $A>B$ holds but
\be
 \exp_\phi(A)-\exp_\phi(B)\le A-B
\nnee
does {\em not} hold.

\vskip 0.1cm

\paragraph{Corollary}
The function $g(u)=\log\exp_\phi(u)$ is increasing and concave,
but {\em not} operator-monotone.

\vskip 0.1cm

\section*{}

\end{document}